\documentclass[12pt]{article}
\nonstopmode
\RequirePackage[colorlinks,citecolor=blue,urlcolor=blue,linkcolor=blue]{hyperref}
\hypersetup{
colorlinks = true,
citecolor=blue,
urlcolor=blue,
linkcolor=blue,
pdfauthor = { Alexey Kuznetsov},
pdfkeywords = { Vandermonde matrix, LU decomposition, Toeplitz matrix,  q-Binomial Theorem, q-Pochhammer symbol,discrete Fourier transform},
pdftitle = {Using q-calculus to study LDLt factorization of a certain Vandermonde matrix},
pdfpagemode = UseNone
}
\usepackage{graphicx,xspace,colortbl}
\usepackage{amsmath,amsthm,amsfonts}
\usepackage{color}
\usepackage{fancybox}
\usepackage{epsfig}
\usepackage{subfig}
\usepackage{pdfsync}
    \def\qed{\hfill$\sqcap\kern-8.0pt\hbox{$\sqcup$}$\\}
    \def\beq{\begin{eqnarray}}
    \def\eeq{\end{eqnarray}}
    \def\beqq{\begin{eqnarray*}}
    \def\eeqq{\end{eqnarray*}}

    \def\c{{\mathbb C}}

    \def\i{{\textnormal i}}

\newtheorem{theorem}{Theorem}

\theoremstyle{definition}

\newtheorem{remark}{Remark}

\title{Using $q$-calculus to study $LDL^T$ factorization of a certain Vandermonde matrix}
\author{
{Alexey Kuznetsov
\footnote{Dept. of Mathematics and Statistics,  York University,
4700 Keele Street, Toronto, ON, M3J 1P3, Canada.  \newline
E-mail:  kuznetsov@mathstat.yorku.ca  \newline Research supported by the Natural Sciences and Engineering Research Council of Canada. } 
 }}
 
 \date{\today}

\begin{document}
\maketitle

\begin{abstract} 
We use tools from $q$-calculus to study $LDL^T$ decomposition of the Vandermonde matrix
$V_q$ with coefficients $v_{i,j}=q^{ij}$. We prove that the matrix $L$ is given as a product of diagonal matrices and the lower triangular Toeplitz matrix $T_q$ with coefficients $t_{i,j}=1/(q;q)_{i-j}$, where $(z;q)_k$ is the q-Pochhammer symbol. 
We investigate some properties of the matrix $T_q$, in particular,  we compute explicitly the inverse of this matrix. 
\end{abstract}

{\vskip 0.15cm}

\section{Introduction and main results}

Let us consider a Vandermonde matrix 
$$
V_q:=
\begin{bmatrix}
    1   & 1 & 1 & 1 & \dots  \\
    1   & q & q^2 & q^3 & \dots  \\
    1  & q^2 & q^4 & q^6 & \dots  \\
    1 &  q^3  & q^6  & q^9 &  \dots  \\
    \vdots & \vdots & \vdots & \vdots & \ddots 
\end{bmatrix}
$$
of size $n\times n$. In the special case when $q=e^{-2\pi \i/n}$, this matrix is called {\it the discrete Fourier transform matrix}. 
Explicit matrix factorizations of the discrete Fourier transform matrix are very important, since they are often used 
in various versions of the Fast Fourier Transform algorithm \cite{VanLoan}. Motivated by this connection, in this note we plan to study the $LDL^T$ factorization of the matrix $V_q$ and to investigate the properties of the factors appearing in such decomposition. The tools and techniques, which are used to prove our results, come from $q$-calculus.

First, let us present several definitions and notation. In what follows, we assume that $n\in {\mathbb N}$ and $q\in \c$. We define the $q$-Pochhammer symbol 
\begin{equation}\label{def_q_Poch}
(z;q)_n:=(1-z)(1-zq)\cdots (1-zq^{n-1}), \;\;\; n\ge 1,
\end{equation}
and $(z;q)_0:=1$. We will  denote by $I$ the $n\times n$ identity matrix. The following matrices of size $n\times n$ will be used frequently in this paper: a lower-triangular Toeplitz matrix $T_q=\{t_{i,j}\}_{0\le i,j\le n-1}$ defined by $t_{i,j}=1/(q;q)_{i-j}$ if $i\ge j$, 
and a diagonal matrix $P_q$ having coefficients $\{(q;q)_i\}_{0\le i \le n-1}$ on the main diagonal, or, more explicitly,
$$
T_q:=
\begin{bmatrix}
    1       & 0 & 0 & 0 & \dots  \\
    \frac{1}{(q;q)_1}   & 1 & 0 & 0 & \dots  \\
        \frac{1}{(q;q)_2}  & \frac{1}{(q;q)_1} & 1 & 0 & \dots  \\
       \frac{1}{(q;q)_3} & \frac{1}{(q;q)_2}  & \frac{1}{(q;q)_1} & 1 &  \dots  \\
    \vdots & \vdots & \vdots & \vdots & \ddots 
\end{bmatrix},
\qquad
P_q:=
\begin{bmatrix}
    1   & 0 & 0 & 0 & \dots  \\
    0   & (q;q)_1 & 0 & 0 & \dots  \\
    0  & 0 & (q;q)_2 & 0 & \dots  \\
    0 &  0  & 0  & (q;q)_3 &  \dots  \\
    \vdots & \vdots & \vdots & \vdots & \ddots 
\end{bmatrix}.
$$
Note that the matrices $T_q$ and $T_{q^{-1}}$ are well-defined for all $q \in \c \setminus {\mathcal A}_n$, where the set ${\mathcal A}_n$ is given by
$$
{\mathcal A}_n:=\{q \in \c : q^j=1 {\textnormal{ for some }} j=1,2,\dots,n-1\}.
$$

In our first result we identify explicitly the matrices appearing in the $LDL^T$ factorization of the Vandermonde matrix $V_q$.
\begin{theorem}\label{thm_main}
Assume that $q \in \c \setminus {\mathcal A}_n$. Then 
$V_q=LDL^T$, where $L=P_q T_q (P_q)^{-1}$ and $D$ is a diagonal matrix having coefficients 
$\{(-1)^i q^{i(i-1)/2} (q;q)_i\}_{0\le i \le n-1}$ on the main diagonal. 
\end{theorem}

In section \ref{sec_proofs} we give a very simple proof of Theorem \ref{thm_main} (our proof is based on the $q$-Binomial Theorem).  
Alternatively, one could derive this result starting from formulas (2.4) and (2.5) in the paper 
\cite{Oruc} by Oruc and Phillips, who use symmetric functions to study LU decomposition of general Vandermonde matrices.

In our second result we present some properties of the Toeplitz matrix $T_q$, including an explicit formula for its inverse. First we define the following two matrices of size $n\times n$: 
\begin{align}\label{def_S_Dq}
S:=
\begin{bmatrix}
    0   & 0 & 0 & 0 & \dots  \\
    1   & 0 & 0 & 0 & \dots  \\
    0  & 1 & 0 & 0 & \dots  \\
    0 &  0  & 1  & 0 &  \dots  \\
    \vdots & \vdots & \vdots & \vdots & \ddots 
\end{bmatrix},
\qquad
D_q:=
\begin{bmatrix}
    1   & 0 & 0 & 0 & \dots  \\
    0   & q & 0 & 0 & \dots  \\
    0  & 0 & q^2 & 0 & \dots  \\
    0 &  0  & 0  & q^3 &  \dots  \\
    \vdots & \vdots & \vdots & \vdots & \ddots 
\end{bmatrix}. 
\end{align}

\begin{theorem}\label{thm1} Assume that $q \in \c \setminus {\mathcal A}_n$. 
\begin{itemize}
\item[(i)] $(T_q)^{-1}=T_{q^{-1}}(I-S)=D_{q^{-1}} T_{q^{-1}} D_q$. 
\item[(ii)] For $m\in {\mathbb N}$ we have
\begin{equation}\label{eqn_m_banded}
T_q D_{q^{-m}} T_{q^{-1}} D_{q^m}=I+\sum\limits_{j=1}^{m-1} \frac{(q^{1-m};q)_j}{(q;q)_j} S^j. 
\end{equation}   
\end{itemize}
\end{theorem}

\begin{remark}
Note that the matrix $H:=(I-S)^{-1}$, which appears in item (i), is a lower triangular Toeplitz matrix having coefficients 
$h_{i,j}=1$ if $i \ge j$ and $h_{i,j}=0$ otherwise. Similarly, the matrix in the right-hand side of \eqref{eqn_m_banded} is a lower-triangular Toeplitz matrix, having $m$ non-zero diagonals: this matrix has coefficient $1$ on the main diagonal and the coefficient $(q^{1-m};q)_j/(q;q)_j$ on the sub-diagonal number $j$, for $1\le j \le m-1$.    
\end{remark}

\section{Proofs}\label{sec_proofs}

The only tool that will be needed for proving Theorems \ref{thm_main} and \ref{thm1} is the q-Binomial Theorem (see \cite{Andrews}[Theorem 10.2.1]), which states that
\begin{equation}\label{eqn_q_binomial_series}
\frac{(az;q)_{\infty}}{(z;q)_{\infty}}=\sum\limits_{j\ge 0} \frac{(a;q)_j}{(q;q)_j} z^j , \;\;\; |q|<1, \; |z|<1. 
\end{equation}
Here $(z;q)_{\infty}:=\prod_{l\ge 0}(1-zq^l)$ and it is clear that this infinite product converges for all $z\in \c$ and 
$|q|<1$. We also record here the following two corollaries of the q-Binomial Theorem, which will be needed later:

\begin{align}
\label{eqn_q_series1}
\frac{1}{(z;q)_{\infty}}&=\sum\limits_{j\ge 0} \frac{z^j}{(q;q)_j}, \;\;\; |q|<1, \; |z|<1, \\
\label{eqn_q_series2}
(z;q)_{\infty}&=\sum\limits_{j\ge 0} \frac{(-1)^j q^{j(j-1)/2}}{(q;q)_j} z^j, \;\;\; |q|<1, \; z\in \c. 
\end{align}

\noindent
{\bf Proof of Theorem \ref{thm_main}:}
First we check that the matrix $L=P_q T_q (P_q)^{-1}$ is a lower-triangular matrix with coefficients 
$$
l_{i,j}=\frac{(q;q)_i}{(q;q)_j (q;q)_{i-j}}, \;\;\; i\ge j.  
$$
Considering an element $(i,j)$ of the matrix $LDL^T$ we see that formula $V_q=LDL^T$ is equivalent to the following  identity: 
for any integers $i,j \ge 0$ 
\begin{equation}\label{finite_sum_identity}
q^{ij}=\sum\limits_{k=0}^{\min(i,j)} \frac{(-1)^k q^{k(k-1)/2}(q;q)_i (q;q)_j}{(q;q)_k (q;q)_{i-k} (q;q)_{j-k}}. 
\end{equation}
We will prove the above identity by writing the Taylor series of the function 
$$
g(u,v):=\frac{(uv;q)_{\infty}}{(u;q)_{\infty} (v;q)_{\infty}}, \;\;\; |u|<1, \; |v|<1, \; |q|<1,
$$
in two different ways. 
First of all, from formula \eqref{eqn_q_binomial_series} we obtain
$$
g(u,v)=\frac{1}{(v;q)_{\infty}} \times \frac{(uv;q)_{\infty}}{(u;q)_{\infty}}=\frac{1}{(v;q)_{\infty}} \sum\limits_{i \ge 0} \frac{(v;q)_i}{(q;q)_i} u^i. 
$$
Using the fact that $(v;q)_i/(v;q)_{\infty}=1/(q^i v;q)_{\infty}$ and expanding this expression in Taylor series in $v$ via \eqref{eqn_q_series1} we conclude that 
\begin{equation}\label{g_series1}
g(u,v)=\sum\limits_{i\ge 0} \sum\limits_{j\ge 0} \frac{q^{ij} u^i v^j }{(q;q)_i (q;q)_j}. 
\end{equation}

On the other hand, we can obtain the series expansion of $g(u,v)$ by applying formulas 
\eqref{eqn_q_series1} and \eqref{eqn_q_series2} in the form
\begin{align*}
(uv;q)_{\infty}&=\sum\limits_{k \ge 0}\frac{(-1)^k q^{k(k-1)/2}}{(q;q)_k} u^k v^k, \\
\frac{1}{(u;q)_{\infty}}&=\sum\limits_{l \ge 0} \frac{u^l}{(q;q)_l}, \\
\frac{1}{(v;q)_{\infty}}&=\sum\limits_{m \ge 0} \frac{v^m}{(q;q)_m}.
\end{align*}
We multiply the above three series expansions and obtain a Taylor series representation in the form 
\begin{equation}\label{g_series2}
g(u,v)=\sum\limits_{k\ge 0} \sum\limits_{l\ge 0} \sum\limits_{m\ge 0}
\frac{(-1)^k q^{k(k-1)/2}}{(q;q)_k (q;q)_l (q;q)_m} u^{k+l} v^{k+m}. 
\end{equation}
Comparing the coefficients in front of the term $u^i v^j$ in both formulas \eqref{g_series1} and \eqref{g_series2} gives us the 
desired result \eqref{finite_sum_identity}.
\qed

\vspace{0.25cm}
\noindent
{\bf Proof of Theorem \ref{thm1}:}
Let us prove the identity $T_q T_{q^{-1}}=(I-S)^{-1}$, which is equivalent to the first equality in item (i) (the second equality 
in (i) follows from formula \eqref{eqn_m_banded} with $m=1$). 
The main idea of the proof is that the Toeplitz  matrix $T_q$ can be expressed in the following form 
\begin{equation}\label{eqn_T_q}
T_q=I+\sum\limits_{j\ge 1} \frac{S^j}{(q;q)_j}.,
\end{equation}
where the matrix $S$ is defined in \eqref{def_S_Dq}. 
The above formula is easy to derive, given that for $1\le j \le n-1$ the coefficients of the matrix $S^j$ have value $1$ on the sub-diagonal number $j$ and value zero everywhere else. In particular, $S^j$ is a zero matrix for $j\ge n$, thus the series in \eqref{eqn_T_q} terminates at $j=n-1$. Similarly, using the identity 
\begin{equation}\label{formula_q_poch_1/q}
(1/q;1/q)_j=(-1)^j q^{-j(j+1)/2} (q;q)_j,
\end{equation}
 and formula \eqref{eqn_T_q} we obtain
\begin{equation}\label{eqn_T_1/q}
T_{q^{-1}}=I+\sum\limits_{j\ge 1} \frac{(-1)^j q^{j(j-1)/2}}{(q;q)_j} (qS)^j. 
\end{equation}
Now, assume that $|q|<1$. Then formulas \eqref{eqn_q_series1} and \eqref{eqn_T_q} give us
\begin{equation}\label{eqn_Tq_product}
T_q=\left[(S;q)_{\infty}\right]^{-1}=(I-S)^{-1} \times (I-qS)^{-1} \times (I-q^2 S)^{-1} \times \cdots.
\end{equation}
Similarly, formulas \eqref{eqn_q_series2} and \eqref{eqn_T_1/q} give us
\begin{equation}\label{eqn_Tq_product2}
T_{q^{-1}}=(qS;q)_{\infty}=(I-qS) \times (I-q^2 S) \times (I-q^3 S) \times \cdots. 
\end{equation}
From the above two identities we see that all the terms $(I-q^i S)$ in the product $T_q T_{q^{-1}}$ are cancelled, except for the first term $(I-S)^{-1}$, thus we obtain $T_qT_{q^{-1}}=(I-S)^{-1}$ for $|q|<1$. We extend this result from $|q|<1$ to 
the general case $q \in \c \setminus {\mathcal A}_n$ by analytical continuation in $q$. 

The proof of formula \eqref{eqn_m_banded} uses the same ideas. Again, first we assume that $|q|<1$. From \eqref{formula_q_poch_1/q} 
  we check that $D_{q^{-m}}  T_{q^{-1}}  D_{q^m}$ is a Toeplitz matrix of the form
\begin{equation*}
D_{q^{-m}} T_{q^{-1}} D_{q^m}=I + \sum\limits_{j\ge 1}  \frac{(-1)^j q^{j(j-1)/2}}{(q;q)_j} (q^{1-m} S)^j
=(q^{1-m}S;q)_{\infty}. 
\end{equation*}
Using the above result and formula \eqref{eqn_Tq_product} we obtain
\begin{equation*}
T_q D_{q^{-m}} T_{q^{-1}} D_{q^m}=\left[(S;q)_{\infty}\right]^{-1} \times (q^{1-m}S;q)_{\infty}=(q^{1-m}S;q)_{m-1}.
\end{equation*}
The desired desired result \eqref{eqn_m_banded} follows by applying \eqref{eqn_q_binomial_series} and analytical continuation in $q$. 
\qed

%

\end{document}